\documentclass[aap]{imsart}
\arxiv{1603.05849}
\RequirePackage[OT1]{fontenc}
\RequirePackage{amsthm,amsmath}
\RequirePackage[numbers]{natbib}
\RequirePackage[colorlinks,citecolor=blue,urlcolor=blue]{hyperref}
\RequirePackage[utf8]{inputenc} % set input encoding (not needed with XeLaTeX)
\RequirePackage{bbm}
\RequirePackage{relsize}
\RequirePackage{
	amsfonts,
	appendix,
	color,
	amssymb,
	graphicx,
	bm,
	stmaryrd,
	mathabx,
	amssymb,
	mathdots,
	esint}
\RequirePackage{amscd,mathrsfs,pifont}
\RequirePackage{eucal,color,empheq}
\RequirePackage{extarrows}
\startlocaldefs
\numberwithin{equation}{section}
\theoremstyle{plain}
\newtheorem{thm}{Theorem}[section]
\newtheorem{cor}{Corollary}[section]
\newtheorem{lemma}{Lemma}[section]

\newtheorem{remark}{Remark}[section]
\newcommand{\bea}{\begin{eqnarray}}
\newcommand{\eea}{\end{eqnarray}}
\newcommand{\R}{\mathbb{R}}

\newcommand{\tr}{\mbox{Tr}}
\newcommand{\pr}{\mathbbm{P}}
\newcommand{\lm}{\lambda_{max}}

\DeclareMathOperator{\erfc}{erfc}
\endlocaldefs
\begin{document}
\begin{frontmatter}
%%%%%%%%%%%%%%%%%%%%%%%%%%%%%%%%%%%%
\title{On the distribution of the largest real eigenvalue for the real
Ginibre ensemble.}
\runtitle{The largest real eigenvalue of a real
Ginibre matrix.}
%%%%%%%%%%%%%%%%%%%%%%%%%%%%%%%%%%%%
\begin{aug}
\author{\fnms{Mihail} \snm{Poplavskyi}\ead[label=e2]{m.poplavskyi@warwick.ac.uk}},
\author{\fnms{Roger} \snm{Tribe}\ead[label=e4]{r.p.tribe@warwick.ac.uk}}
\and
\author{\fnms{Oleg} \snm{Zaboronski}\corref{}
\ead[label=e5]{olegz@maths.warwick.ac.uk}}

\address{MP, RT, OZ: Department of Mathematics, University of Warwick, Coventry CV4~7AL, UK\\
\printead{e2}\\
\phantom{E-mail:\ }\printead*{e4}\\
\phantom{E-mail:\ }\printead*{e5}
}

\end{aug}

%%%%%%%%%%%%%%%%%%%%%%%%%%%%%%%%%%
%ABSTRACT
%%%%%%%%%%%%%%%%%%%%%%%%%%%%%%%%%%
\begin{abstract}
Let $\sqrt{N}+\lambda_{max}$ be the largest real eigenvalue of a random $N\times N$ matrix
with independent $N(0,1)$ entries (the `real Ginibre matrix'). We study the large deviations
behaviour of the limiting $N\rightarrow \infty$ distribution $\pr[\lm<t]$ of the shifted maximal real eigenvalue $\lambda_{max}$. 
In particular, we prove that the right tail of this distribution is Gaussian: for $t>0$, 
\[
\pr[\lm<t]=1-\frac{1}{4}\erfc(t)+O\left(e^{-2t^2}\right).
\]
This is a rigorous confirmation of the corresponding result of \cite{forrester_nagao}. We also prove that the left tail is exponential, with correct asymptotics up to $O(1)$: for $t<0$,
\[
\pr[\lm<t]=
e^{\frac{1}{2\sqrt{2\pi}}\zeta\left(\frac{3}{2}\right)t+O(1)},
\] 
where $\zeta$ is the Riemann zeta-function. 

Our results have implications for interacting particle systems.
The edge scaling limit of the law of real eigenvalues for the real Ginibre ensemble is a rescaling of a fixed time distribution
of annihilating Brownian motions (ABM's) with the step initial condition, see \cite{bg_mp_rt_oz}.
Therefore, the tail behaviour of the distribution of $X_s^{(max)}$ - the position of the rightmost annihilating particle at fixed time $s>0$ - can be
read off from the corresponding answers for
$\lm$ using $X_s^{(max)}\stackrel{D}{=} \sqrt{4s}\lambda_{max}$.
\end{abstract}
%%%%%%%%%%%%%%%%%%%%%%%%%%%%%%%%%%%%%%%%%%%%%%%%%%%%%
%%%%%%%%%%%%%%%%%%%%%%%%%%%%%%%%%%%%%%%%%%%%%%%%%%%%%

%%%%%%%%%%%%%%%%%%%%%%%%%%%%%%%%%%%
%%%%%%%%%%%%%%%%%%%%%%%%%%%%%%%%%%%
\end{frontmatter}
%%%%%%%%%%%%%%%%%%%%%%%%%%%%%%%%%%%
%BODY OF THE PAPER
%%%%%%%%%%%%%%%%%%%%%%%%%%%%%%%%%%%
\section{Introduction and the main result}\label{sec1}
%%%%%%%%%%%%%%%%%%%%%%%%%%%%%%%%%%%
The laws describing the distribution of eigenvalues of large self-adjoint random matrices near
the spectral edge exhibit a large degree of universality. One of the examples is the celebrated family of Tracy-Widom distributions for the largest eigenvalue
$F_{\beta}, \beta=1,2,4$, which has been originally discovered in the context of Gaussian orthogonal, unitary
and symplectic ensembles \cite{tw1}, \cite{tw2}. These distributions also describe the scaling limit of the distribution of the largest eigenvalue for non-Gaussian invariant ensembles \cite{deift}, \cite{goev}
as well as non-invariant ensembles of random matrices with independent non-Gaussian entries \cite{sosh}.
Moreover, $F_\beta$ appears in the large number of statistical models not directly related to random matrices, such as
random permutations, growth models belonging KPZ universality class and related asymmetric exclusion processes, see \cite{tw_review}
for a recent review.

There is a growing body of evidence that the extreme statistics of eigenvalues for the real Ginibre ensemble \cite{ginibre} also give rise to a
new universality class which is relevant beyond random matrix theory.
Recall that the real Ginibre ensemble is a Gaussian measure on the set
of random real $N\times N$ matrices such that the matrix elements are
independent $N(0,1)$ random variables. The law of the real eigenvalues
in the $N\rightarrow \infty$ edge scaling limit is governed by a Pfaffian point process discovered independently in \cite{forrester_nagao},
\cite{borodin_sinclair} and \cite{sommers}. Within random matrix theory, this Pfaffian point process turns out
to be universal, in the sense that it holds for a large class of non-Gaussian
ensembles of real non-symmetric matrices, see \cite{tao_vu}, Corollary $15$.\footnote{The cited universality statement
 covers both the edge and the bulk scaling limits, but only the former is needed for the current
investigation.} The universality of the extreme eigenvalue statistics can then be shown to follow from the universality
of the local correlations functions via the Fredholm Pfaffian representation, see \cite{rider_sinclair_2014}. 
Outside random matrix theory,
the same Pfaffian point process describes the law of instantaneously annihilating Brownian
motions on the line started with half-space initial conditions, see \cite{bg_mp_rt_oz}. 
\footnote{The bulk scaling limit of the law of the real eigenvalues for the real Ginibre
ensemble coincides with the fixed time law for the annihilating Brownian motions
corresponding to the maximal entrance law, \cite{roger_oleg}.}
Therefore, the scaling limit of the law of the
largest real eigenvalue for the real Ginibre ensemble is the same as the law of the rightmost annihilating particle at a fixed 
time. In this paper we will analyse the tails of this distribution and show that
\bea\label{eqn:1}
\pr(\lambda_{max}<t) \xlongequal{t\rightarrow -\infty} e^{\frac{\zeta(3/2)}{2\sqrt{2\pi}}t+O(1)},
\eea
where $\lambda_{max}$ is the position of the maximal eigenvalue measured from the right edge
of the spectrum for the real Ginibre ensemble or it is the scaled position $X^{(max)}_s/\sqrt{4s}$ of the rightmost particle for annihilating Brownian motions
at time $s$. 

An instance of edge statistics (\ref{eqn:1}) outside random matrix theory $and$ Pfaffian point processes is for the symmetric exclusion process on $\mathbb{Z}$ with half-filled initial conditions. That is, all sites
to the left of zero are occupied at time $s=0$, all sites to the right of zero are empty;
at time $s>0$, particles hop onto empty neighboring sites at rate $1$. Let $R_s$ be the position of the
rightmost particle at time $s$. It is shown in \cite{pavel1} that
\bea
\pr(R_s=0) \xlongequal{s\rightarrow \infty} e^{-\frac{\zeta(3/2)}{\sqrt{\pi}}\sqrt{s}+o(\sqrt{s})},
\eea
which coincides with (\ref{eqn:1}) if we identify distance $|t|$ with the diffusive scale $\sqrt{8s}$.
To verify that this agreement is not accidental, one must compare $\pr(\lambda_{max}<t)$
and $\pr(R_{t^2/8}=0)$ for all values of $t$, not just in the large-$|t|$ limit, but the very possibility that there is a relation between the real Ginibre ensemble and the symmetric exclusion process is intriguing.

Let us stress that formula (\ref{eqn:1}) is easy to guess given results on bulk gap probabilities 
described in Forrester \cite{forrester}. Building on the paper \cite{derrida} on gap probability
for annihilating Brownian motions and using the relation between the bulk statistics for the real Ginibre ensemble
and annihilating Brownian motions with the maximal entrance law initial conditions \cite{roger_oleg},
Forrester argues that the bulk scaling limit of the gap probability for the real Ginibre
ensemble satisfies (\ref{eqn:1}) if $|t|$ is interpreted as the size of the gap. As the width of the transition region
near the edge of the spectrum is of order one for the real Ginibre ensemble, it is only natural to guess
that (\ref{eqn:1}) stays valid at the edge with edge effects showing only as $O(1)$ terms.
The aim of the current paper is to verify these heuristics rigorously.

The starting point of our investigation is the following Fredholm determinant representation
of $\pr[\lm<\alpha]$ due to Rider and Sinclair \cite{rider_sinclair_2014}\footnote{The paper
\cite{rider_sinclair_2014} contains some minor computational errors, which led to an incorrect expression for
$\pr(\lm<\alpha)$. Here we present the corrected version. An interested reader is referred to the appendix for
a discussion of the errors and their correction.}:
\begin{thm}[B. Rider and C. D. Sinclair, 2014]\label{thm:rs}
Introduce the integral operator $T$ with kernel
\bea
T(x,y)=\frac{1}{\pi} \int_0^\infty e^{-(x+u)^2} e^{-(y+u)^2}du.
\eea
Let $\chi_t$ be the indicator of $(t,\infty)$. Then, 
\bea 
\pr[\lm<t]=\sqrt{\det(I-T\chi_t)\Gamma_t},
\eea
where $\Gamma_t$ is defined as follows. Set $g(x)=\frac{1}{\sqrt{\pi}}e^{-x^2}$ 
, $G(x)=\int_{-\infty}^x g(y)dy$
and denote by $R(\cdot,\cdot)$ the kernel of the operator $(I-T\chi_t)^{-1}-I$. Then
\begin{eqnarray*}
\Gamma_t&=&\left(1-\frac{1}{2}a_t\right) \left(1-\frac{1}{2}\int_{-\infty}^t R(x,t+)dx\right)\\
&&+\frac{1}{4}\left(1-b_t\right)\left(\int_{-\infty}^t (I-T\chi_t)^{-1} g(x)dx-1\right)
\end{eqnarray*}
for $a_t=\int_t^\infty G(x) (I-T\chi_t)^{-1} g(x)dx$ and $b_t=(I-T\chi_t)^{-1}G(t)$.
\end{thm}
The reader is referred to the original paper for precise definitions of the quantities entering the theorem.
We mention briefly only that $\det(I-T\chi_t)$ should be understood as the Fredholm determinant
of the operator $T$ acting on $L^{2}(t,\infty)$ and that $T$ is the square of a Hilbert-Schmidt operator, 
which implies that $T$ is positive and trace class.  

The apparently complicated expression for the factor $\Gamma_t$ in the statement of the Rider-Sinclair theorem 
is actually rather simple, and cancellations occur that imply that $\Gamma_t=(1-a_t)$. Therefore 
\begin{lemma}\label{thm:rsmod}
\bea \label{eq:rsmod}
\pr(\lm<t)=\sqrt{\det(I-T\chi_t)(1-a_t)}, 
\eea
where the operator $T$ and the function $a_t$ are defined in the statement of Theorem \ref{thm:rs} 
\end{lemma}
In what follows we are going to use the modified version (\ref{eq:rsmod}) of the Rider-Sinclair result.

As it turns out, the asymptotic analysis of (\ref{eq:rsmod}) can be greatly simplified by assigning a probabilistic interpretation
to all the factors. In particular, the interpretation allowed us to control
the asymptotic expansion beyond the leading order, which is hard to do using
purely operator theoretic methods.
Let $\left(B_n,n\geq 0\right)$ be the discrete time random  walk with Gaussian $N(0,1/2)$ increments started at zero.
Let
\begin{equation} \label{eq:tau0}
\left\{ \begin{array}{rcl}
\tau_t &=&\inf_{n> 0} \{2n-1:~B_{2n-1}\leq t\},\\
\tau_0&=&\inf_{n> 0} \{2n:~B_{2n} \geq 0\}.
\end{array} \right. 
\end{equation}
In words: $\tau_0$ is the smallest $even$ time such that $B_{\tau_{0}}\geq 0$, $\tau_t$ is the smallest $odd$ time such that
$B_{\tau_t}\leq t$. 
Also, let 
\bea
I_{\tau_0}=\inf\{B_s:s~is~odd,~s\leq \tau_0 \}
\eea
be the infimum of the random walk $(B_s)_{\geq 0}$ taken over all odd times not exceeding the exit time $\tau_0$.
Then we have the following probabilistic restatement of Theorem \ref{thm:rs}:
\begin{thm} \label{thm:prob}
\bea\label{eq:prob}
\\
\nonumber
\pr(\lm<t)=\sqrt{\pr(\tau_t<\tau_0)}
e^{-\frac{1}{2}\mathbb{E}\left(\left(I_{\tau_0}-t\right)_{+}\delta_0\left(B_{\tau_0}\right)\right)},
\eea
where $x_{+}:=\max(x,0)$ is the positive part of a real number. 
\end{thm}
\noindent
{\bf Remark.} We use the expression $\mathbb{E}(X\delta_{y}(Y))$ to mean a continuous Lebesgue density for the measure 
$\mathbb{E}(X\mathbbm{1}(Y\in dy))$ evaluated at $y$. If $y=0$ we sometimes rewrite $\mathbb{E}(X\delta_{0}(Y))$ as 
$\mathbb{E}(X\mathbbm{1}(Y\in d0))$.\\
\\
Note that (\ref{eq:prob}) can be rewritten in a more useful way when $t<0$ as
\bea\label{eq:probm}
\\
\nonumber
\pr(\lm<t)=\sqrt{\pr(\tau_t<\tau_0)} \;
e^{\frac{t}{2} \; \mathbb{E}\left(\delta_0\left(B_{\tau_0}\right)\right)}
e^{- \frac{1}{2} \mathbb{E}\left(\max(t,I_{\tau_0})\delta_0\left(B_{\tau_0}\right)\right)}.
\eea
Our main result below will follow from the above probabilistic representation via an application of
some general results for random walks. In particular it can be easily shown that $\mathbb{E}\left(\delta_0\left(B_{\tau_0}\right)\right)
= \zeta(3/2)/\sqrt{2\pi}$ (see Lemma \ref{keyconstant}), which then controls the first term in the asymptotics for $t<0$. 
The following theorem details the asymptotics up to $O(1)$. 
\begin{thm}\label{thm:asympt}
For $t>0$,
\bea\label{eq:tpos}
\pr(\lm<t)=1-\frac{1}{4}\erfc(t)+O(e^{-2t^2}).
\eea
For $t<0$,
\bea\label{eq:tneg}
\pr(\lm<t)=\exp\left(-\frac{\zeta(3/2)}{2\sqrt{2\pi}}|t|+O(1)\right).
\eea
\end{thm}
For $t<0$, the above statement provides a rigorous justification of (\ref{eqn:1}). The right tail
asymptotic (\ref{eq:tpos}) is actually well known \cite{forrester_nagao} since the probability of finding an eigenvalue very
far to the right of the spectral edge is approximately equal to the level density. This part of Theorem \ref{thm:asympt}
should be considered as a test of the Rider-Sinclair answer complemented by careful bounding of error terms.

As has been already mentioned, the edge scaling limit of the law of real eigenvalues for the real Ginibre
ensemble coincides up to a Brownian rescaling with the law of annihilating Brownian motions on the real line started with  the step
initial condition \cite{bg_mp_rt_oz}. The step initial condition corresponds to the maximal entrance law restricted to $x<0$ and zero density
of particles for $x>0$. The maximal entrance law can be constructed as the limit of the homogeneous Poisson
point process initial conditions with intensity diverging to infinity, see
  \cite{bg_mp_rt_oz} for precise definitions and details. Therefore, Theorem \ref{thm:asympt}
yields a simple corollary:
\begin{cor}\label{cor:abm}
Consider the system of instantaneously annihilating Brownian motions on the real line started with the step initial condition. 
Let $X_s^{(max)}$ be the position of the rightmost particle at a fixed time $s>0$. Then
\[
\pr(X_s^{(max)}<x)=1-\frac{1}{4}\erfc\left(\frac{x}{\sqrt{4s}}\right)
+O\left(e^{-\frac{x^2}{2s}}\right)
\]
for $x/\sqrt{s}\rightarrow \infty$, while for $x/\sqrt{s} \to -\infty$, 
\[
\pr(X_s^{(max)}<x)=
e^{\frac{1}{2\sqrt{2\pi}}\zeta\left(\frac{3}{2}\right)\frac{x}{\sqrt{4s}}+O(1)}.
\]
\end{cor}
The above result complements the recent study \cite{pavel2} of annihilating Brownian motions near the edge of the
distribution, where the average number of particles in the positive half-line $(0,\infty)$ has been calculated.

The correspondence between the real eigenvalues and annihilating particles can also be used
for an intuitive explanation of  Theorem \ref{thm:asympt}: the right tail of $\pr(\lm<t)$ is Gaussian 
as it corresponds to the probability that a Brownian particle travels distance $t$ during the time interval $[0,1/4]$. The left tail cannot be
thinner than exponential, as the probability of the event $\{\lm<t\}$ for a negative $t$ of large magnitude can be bounded
below by the intersection of the following $O(|t|)$ independent events: particles stay within each of $|t|$
boxes of size $1$ and completely annihilate each other by the time $s=1/4$; particles with initial positions to the left of $t$ do not enter the interval
$(t,\infty)$ before time $1/4$. Unfortunately, we could not extract the exact rate of the exponential decay from arguments like this.

The rest of the paper is organised as follows. In Section \ref{sec2} we will prove Lemma \ref{thm:rsmod}
and Theorem \ref{thm:prob}. The asymptotic analysis of $\pr(\lambda_{max}<s)$ is carried out in Section \ref{sec3}.
All technical lemmas needed for the proof of our main result are proved in Appendix A. In Appendix B
we explain how to fix the minor errors we found in the original statement of the Rider-Sinclair Theorem. 

%%%%%%%%%%%%%%%%%%%%%%%%%%%%%%%%%%%%%%%%%%%%%%%%
%%%%%%%%%%%%%%%%%%%%%%%%%%%%%%%%%%%%%%%%%%%%%%%%

%%%%%%%%%%%%%%%%%%%%%%%%%%%%%
\section{The proof of Lemma \ref{thm:rsmod} and Theorem \ref{thm:prob}}\label{sec2}
%%%%%%%%%%%%%%%%%%%%%%%%%%%%%
Let us examine the $\Gamma$-factor defined in the statement of the Rider-Sinclair Theorem \ref{thm:rs}.
Throughout the paper we will use the following basic property of the integral operator $T$:
\begin{lemma}\label{thm:l1}
For any $t \in \R$, the operator $T$ acting on $L^{2}(t,\infty)$ is positive definite. Its spectral radius $\rho_t$
is bounded away from $1$: for any $t>0$,
\bea\label{eq:pbound}
\rho_t \leq \frac{1}{8}e^{-2t^2}.
\eea
There exists $C_0>0$ such that for any $t<0$,
\bea\label{eq:nbound}
\rho_t\leq \exp{\left[-\frac{C_0}{1+t^2}\right]}.
\eea
\end{lemma}
Therefore, for any $t\in \R$, the resolvent of the operator $T\chi_t$ can be expanded into an absolutely
convergent power series,
\[
(I-T\chi_t)^{-1}=\sum_{n=0}^\infty (T\chi_t)^n.
\]
Using the explicit definitions of the integral operator $T$ and functions $G,g$ from Theorem \ref{thm:rs},
\bea\label{eq:at}
a_t:=\int_t^\infty G(x) (I-T\chi_t)^{-1} g(x)dx=\sum_{m=0}^\infty p_m,
\eea
where
\bea
p_m=\int_{\R^{2m+2}} \prod_{k=0}^m \frac{dy_kdu_k}{\pi} e^{-y_m^2-\sum_{n=1}^m\left((u_n-y_n)^2+(y_{n-1}-u_n)^2\right)-(y_0-u_0)^2}
\nonumber \\
\times \left(\prod_{p=0}^m\chi(y_p\geq t)
\prod_{q=1}^m\chi(u_q\leq  0) \chi(u_0\geq 0)\right),
\label{eq:pm}
\eea
where $\chi(\cdot \geq a)$, $\chi(\cdot \leq a)$ are the characteristic functions of $[a,\infty)$ and $(-\infty,a]$
correspondingly.
The expression for $p_m$ quoted above arises from some simple changes of variables designed to 
re-write in terms of the density of a segment
of the Gaussian random walk introduced in Section \ref{sec1}. The integrand in (\ref{eq:pm}) can be
identified with the density for the initial
segment of the walk of length $2m+2$
$
(B_1,B_2,\ldots, B_{2m+2})$ at the point
$(y_m,u_m,y_{m-1},u_{m-1},\ldots y_{0},u_{0})$.
This reveals that
\[
p_m=\pr\left[B_{2k+1}\geq t, k=0,1,\ldots , m; B_{2l}\leq 0,l=1,2,\ldots m;B_{2m+2}\geq 0\right].
\]
Equivalently, in terms of exit times $\tau_t, \tau_0$ defined in (\ref{eq:tau0}),
\[
p_m=\pr\left[\tau_0=2m+2; \tau_t>\tau_0 \right],~m\geq 0.
\]
Substituting this formula into (\ref{eq:at}) and summing over $m$'s, we find that
\bea\label{eq:aprob}
a_t=\pr[\tau_t>\tau_0].
\eea
Note that the above expression is well defined, as the exit time $\tau:=\tau_0\wedge \tau_t$
is a finite random variable, $\pr[\tau<\infty]=1$.

Next,
\bea\label{eq:bt}
b_t:=(I-T\chi_t)^{-1} G(t)=\sum_{m=0}^\infty q_m,
\eea
where
\bea
q_m=\int_{\R^{2m+1}} \prod_{k=1}^m \frac{dy_kdu_k}{\pi}\frac{du_{m+1}}{\sqrt{\pi}} e^{-u_1^2-\sum_{n=1}^m\left((y_n-u_n)^2+(u_{n+1}-y_n)^2\right)}
\nonumber \\
\times \left(\prod_{p=1}^m\chi(y_p\leq 0)
\chi(u_p\geq  t) \chi(u_{m+1}\leq t)\right),
\label{eq:qm}
\eea
Consider the density for segment $(B_1,B_2,\ldots, B_{2m+1})$ of the random walk at the point $(u_1,y_1,\ldots u_{m},y_{m},u_{2m+1})$.
Re-writing (\ref{eq:bt}) as an expectation over this segment we find that
\[
q_m=\pr\left[B_{2k+1}\geq t, k=0,1,\ldots , m-1;B_{2l}\leq 0,l=1,2,\ldots m;B_{2m+1}\leq t\right].
\]
In terms of exit times this is
\[
q_m=\pr\left[\tau_t=2m+1, \tau_t<\tau_0 \right], ~m\geq 0.
\]
Substituting this expression into (\ref{eq:bt}) and summing over $m$, we find that
\bea\label{eq:bprob}
b_t=\pr[\tau_t<\tau_0]=1-a_t.
\eea
In a very similar fashion we find that
\bea\label{eq:cprob}
\int_{-\infty}^t R(x,t)dx=\pr[\tau_t>\tau_0]=a_t
\eea
and
\bea\label{eq:dprob}
\int_{-\infty}^t (I-T\chi_t)^{-1}g(x)dx=\pr[\tau_t<\tau_0]=1-a_t.
\eea
Substituting (\ref{eq:aprob}), (\ref{eq:bprob}), (\ref{eq:cprob}) and (\ref{eq:dprob}) into the expression for $\Gamma_t$
presented in Theorem \ref{thm:rs} gives
\[
\Gamma_t=(1-a_t/2)^2-a_t^2/4=1-a_t=\pr[\tau_t<\tau_0]. 
\]
Lemma \ref{thm:rsmod} and the pre-factor in formula (\ref{eq:prob}) of Theorem \ref{thm:prob} are verified.

To check the exponent in the right hand side of (\ref{eq:prob}), we use Lemma \ref{thm:l1} to
justify the application of the trace-log formula to the Fredholm determinant appearing in the Rider-Sinclair Theorem
and the subsequent expansion of the logarithm:
\[
\log \det(I-T\chi_t)=\tr \log(I-T\chi_t)=-\sum_{m=1}^{\infty} \frac{1}{m} \tr (T\chi_t)^m.
\]
Explicitly,
\bea\label{eq:log}
\log\det (I-T\chi_t)=-\sum_{m=1}^{\infty} \frac{1}{m}r_m(t),
\eea
where
\bea
r_m(t)=\int_{\R^{2m}} \prod_{k=1}^m \frac{dy_kdu_k}{\pi} e^{-\sum_{n=1}^m\left((u_n-y_n)^2+(y_{n+1}-u_n)^2\right)}
\nonumber \\
\times \left(\prod_{p=1}^m\chi(u_p\leq 0)
\chi(y_p\geq  t)\right),
\label{eq:rm}
\eea
where $y_{m+1}\equiv y_1$. Differentiating $r_m$ with respect to $t$ gives a sum of $m$
terms which are all identical due to the cyclic symmetry of the integrand.
Interpreting each of the terms using the Gaussian random walk of length $2m$, we
find that
\[
\frac{dr_m}{dt}(t)=-m\pr[\tau_0=2m; \tau_t>\tau_0;B_{\tau_0}\in d0].
\]
Substituting the above into the $t$-derivative of (\ref{eq:log}) one gets 
\begin{eqnarray*}
\frac{d}{dt}\log \det(I-T\chi_t)&=&\pr[ \tau_t>\tau_0;B_{\tau_0}\in d0]\\
&=&
\mathbb{E}( \mathbbm{1}(I_{\tau_0}>t)\delta_0( B_{\tau_0})),
\end{eqnarray*}
where we used that $\{\tau_t>\tau_0\}=\{I_{\tau_0}>t\}$.
Integrating the above expression over the interval $(t,+\infty)$ and applying the boundary
condition $\det(I-T\chi_t)\mid_{t=+\infty}=1$, we arrive at the claimed expression for the exponent:
\bea
\log \det(I-T\chi_t)=-\mathbb{E}((I_{\tau_0}-t)_{+}\delta_0( B_{\tau_0})).
\eea
Formula (\ref{eq:prob}) is proved.

Finally, equation (\ref{eq:probm}) follows from (\ref{eq:prob}) due to the following elementary identity:
for any $t,~I$,
\[
t+(I-t)_{+}=\max(t,I).
\]

Theorem \ref{thm:prob} is proved.

%%%%%%%%%%%%%%%%%%%%%%%%%%%%%%%%%%%%%%%%%%%%%
\section{The proof of Theorem \ref{thm:asympt}}\label{sec3}
%%%%%%%%%%%%%%%%%%%%%%%%%%%%%%%%%%%%%%%%%%%%%
\subsection{Asymptotic expansion for $t<0$: the proof of (\ref{eq:tneg})}

A natural starting point for our analysis is equation (\ref{eq:probm}). 
Firstly, let us investigate the coefficient of the $O(t)$ term, $\mathbb{E}(\delta_0(B_{\tau_0}))$. 
The event $\{B_{\tau_0}\in d0\}$ depends only on the position of the random walk at even times.
Averaging over positions at odd times we find that 
\[
\pr[B_{\tau_0}\in d0]=\pr[\tilde{B}_{t_0} \in d0],
\]
where $(\tilde{B}_s)_{s \geq 0}$ is a Gaussian random walk with $N(0,1)$-increments
and $t_0=\inf \{s=1,2,\ldots: \tilde{B}_s>0\}$.
Therefore
\begin{eqnarray*}
\mathbb{E}(\delta_0(B_{\tau_0}))&=&\sum_{n=1}^\infty \pr [t_0=n,\tilde{B}_n \in d0]\\
&=&\sum_{n=1}^\infty \pr [\tilde{B}_1<0; \tilde{B}_2<0; \tilde{B}_{n-1}<0;\tilde{B}_n \in d0]\\
&=&\sum_{n=1}^\infty \pr [\tilde{B}_1<0; \tilde{B}_2<0; \tilde{B}_{n-1}<0\mid \tilde{B}_n =0]\pr [ \tilde{B}_n \in d0]
\end{eqnarray*}
The last expression can be evaluated using the following remarkable combinatorial lemma:
\begin{lemma}\label{thm:comb}
\bea
\pr [\tilde{B}_1<0; \tilde{B}_2<0; \tilde{B}_{n-1}<0\mid \tilde{B}_n =0]=\frac{1}{n}.
\eea
\end{lemma}
This lemma follows from a more general combinatorial result concerning the total time a random
walk spends above zero, see \cite{feller}, Chapter XII, Section 6. For the sake of completeness, a very simple proof of Lemma \ref{thm:comb}
is presented in the Appendix. 
We conclude that
\bea \label{keyconstant}
\mathbb{E}(\delta_0(B_{\tau_0}))=\sum_{n=1}^\infty \frac{\pr [ \tilde{B}_n \in d0]}{n}=
\sum_{n=1}^\infty \frac{1}{\sqrt{2\pi n^3}}=\frac{\zeta(3/2)}{\sqrt{2\pi}}.
\eea
The leading order asymptotic for $\log \pr[\lambda_{max}<t]$ is derived.

To finish the derivation of equation (\ref{eq:tneg}) from (\ref{eq:probm}) we have to show that
\bea\label{eq:step}
\log \pr[\tau_t<\tau_0]-\mathbb{E}(\max(t,I_{\tau_0})\delta_0(B_{\tau_0}))=O(1).
\eea
By a Brownian motion analogue, it is clear that each of the terms on the left hand side
of (\ref{eq:step}) is $O(\log|t|)$. The main challenge is to show that the logarithms cancel.
This cancellation follows from the following two results:
\begin{lemma}\label{thm:exit}
There exists a positive constant $C$ such that for a sufficiently large $|t|$,
\end{lemma}
\bea\label{eq:mart}
\frac{1}{\sqrt{2}|t|}\geq \pr[\tau_t<\tau_0]\geq \frac{1}{\sqrt{2}|t|} (1-C|t|^{-1/2}).
\eea
\begin{lemma}\label{thm:lotov}
Consider a Gaussian random walk with $N(0,1)$ increments.
There exists a positive constant $\mu$ such that for $y<0$,
\bea
\mathbb{E}_y(\delta_0(B_{\tau_0}))=\sqrt{2}+O(e^{-|y|\mu}).
\eea
Here the subscript $y$ means that the random walk is started from $y$ and $\tau_0=\inf(n: B_n>0)$.
\end{lemma}
Lemma \ref{thm:exit} can be easily proved using martingale methods, see the Appendix. In contrast, Lemma \ref{thm:lotov}
is a corollary of Theorem 4 of \cite{lotov}, which is a result of an intricate asymptotic analysis of Wiener-Hopf equations associated with exit problems for
Gaussian random walks.

It follows from Lemma \ref{thm:exit} that
\bea\label{eq:step1}
\log \pr[\tau_t<\tau_0]=-\log|t|+O(1).
\eea
The second term on the left hand side of (\ref{eq:step}) can be simplified using integration by parts. The result is
\[
\mathbb{E}(\max(t,I_{\tau_0})\delta_0(B_{\tau_0}))=-\int_{t}^0 dy \mathbb{E}(\mathbbm{1}_{I_{\tau_0}<y}\delta_0(B_{\tau_0}))+
\int_0^\infty dy \mathbb{E}(\mathbbm{1}_{I_{\tau_0}>y}\delta_0(B_{\tau_0})).
\]
The second term on the right hand side of the above identity is $t$-independent and is finite as the the right tail of the distribution of $I_{\tau_0}$
conditioned on $B_{\tau_0}=0$ can be shown to have Gaussian decay. The first term can be evaluated using Lemma \ref{thm:lotov}:
let $L<0$ be a $t$-independent constant. Then
\begin{eqnarray*}
\int_{t}^0 dy \mathbb{E}(\mathbbm{1}_{I_{\tau_0}<y}\delta_0(B_{\tau_0}))&=&
\left(\int_{t}^L+\int_L^0\right) dy\mathbb{E}
(\mathbbm{1}_{I_{\tau_0}<y}\delta_0(B_{\tau_0}))\\
&=&\int_{t}^L dy \mathbb{E}(\mathbbm{1}_{I_{\tau_0}<y}\delta_0(B_{\tau_0}))+O(1)\\
&=&\int_{t}^L dy \mathbb{E}(\mathbbm{1}_{I_{\tau_0}<y}\mathbb{E}_{B_{\tau_y}}(\delta_0(B_{\tau_0})))+O(1)\\
&\stackrel{L. \ref{thm:lotov}}{=\joinrel=\joinrel=}&
\int_{t}^L dy \mathbb{E}(\mathbbm{1}_{I_{\tau_0}<y}(\sqrt{2}+O(e^{-\mu|B_{\tau_y}|})))+O(1)\\
&\stackrel{L. \ref{thm:exit}}{=\joinrel=\joinrel=}&\log(|t|)+O(1).
\end{eqnarray*}
We conclude that
\bea\label{eq:step2}
\mathbb{E}(\max(t,I_{\tau_0})\delta_0(B_{\tau_0}))=-\log(|t|)+O(1).
\eea
Substituting (\ref{eq:step1}) and (\ref{eq:step2}) into the left hand side of (\ref{eq:step}) we observe that the 
logarithms cancel and we are left with terms of order $1$, as claimed. Theorem \ref{thm:asympt} is proved for $t<0$.

%%%%%%%%%%%%%%%%%%%%%%%%%%%%%%%%%%%%%%%%%%%%%
\subsection{Asymptotic expansion for $t>0$: the proof of (\ref{eq:tpos})}
%%%%%%%%%%%%%%%%%%%%%%%%%%%%
The $t>0$ case is best tackled starting directly from the Rider-Sinclair theorem with the pre-factor given by Lemma
\ref{thm:rsmod}. Recall that the operator $T\chi_t$ is a positive trace class operator on $L^2(t,\infty)$
with spectral radius strictly less than 1, see Lemma \ref{thm:l1}. Therefore,
\[
\pr[\lambda_{max}<t]=\sqrt{(1-a_t)}e^{-\frac12 Q_t},
\]
where $Q_t=\sum_{m=1}^\infty \frac{1}{m} \tr (T\chi_t)^m$. An explicit calculation 
(see \ref{trace}) shows that
\[
\tr (T\chi_t)\leq \frac{1}{8}e^{-2t^2}.
\]
This leads to the following estimate:
\begin{eqnarray*}
0\leq Q_t&=&\sum_{m=1}^\infty \frac{1}{m} \tr (T\chi_t)^m
\leq \sum_{m=1}^\infty \frac{1}{m} (\tr (T\chi_t))^m \\
&=&-\log(1-\tr(T\chi_t))\leq -\log\left(1-\frac{1}{8} e^{-2t^2}\right).
\end{eqnarray*}
Therefore we conclude that
\bea\label{eq:prxp}
\pr[\lambda_{max}<t]=\sqrt{(1-a_t)}(1+O(e^{-2t^2})).
\eea
To calculate the pre-factor, notice that
\[
T\chi_t g(y)=g(y)\int_{t}^\infty dx\int_0^\infty\frac{dz}{\pi} e^{-2(x^2+z^2+xz+yz)} 
\leq \frac{1}{8\sqrt{2\pi} t} e^{-2t^2}g(y). 
\]
so,
\[
a_t=\int_t^\infty G(x)\sum_{k=0}^\infty (T\chi_t)^kg(x)dx=\int_t^{\infty} G(x)g(x)dx+\rho_t,
\]
where for suitably large $t$
\begin{eqnarray*}
0\leq \rho_t \leq \int_t^{\infty} G(x)g(x)dx \sum_{m=1}^\infty \left(\frac{1}{8\sqrt{2\pi} t}e^{-2t^2}\right)^m
=\int_t^{\infty} G(x)g(x)dx\cdot
\frac{\frac{1}{8\sqrt{2\pi} t}e^{-2t^2}}{1-\frac{1}{8\sqrt{2\pi} t}e^{-2t^2}}.
\end{eqnarray*}
Also notice that
\[
\int_t^{\infty} G(x)g(x)dx=\int_t^{\infty} G(x)G'(x)dx=\frac{1}{2}(1-G(t)^2)=\frac{\erfc(t)}{2}-\frac{\erfc^2(t)}{8}.
\]
Therefore,
\bea\label{eq:aat}
a_t=\frac{1}{2}\erfc(t)(1+O(e^{-{t^2}})).
\eea
Substituting (\ref{eq:aat}) into (\ref{eq:prxp}) we find that for $t>0$,
\[
\pr[\lambda_{max}<t]=1-\frac{1}{4}\erfc(t)+O(e^{-2t^2}).
\]
Formula (\ref{eq:tpos}) is proved.
%%%%%%%%%%%%%%%%%%%%%%%%%%%%%%%%%%%%%%%%
\section*{Acknowledgements}
We are grateful to Paul Krapivsky and Gernot Akemann for many illuminating discussions. 
This research was supported by EPSRC grant No.\ EP/K011758/1 (M.P. and R.T.) and  a Leverhulme Trust Research Fellowship (O.Z.).
%%%%%%%%%%%%%%%%%%%%%%%%%%%%

%%%%%%%%%%%%%%%%%%%%%%%%%%%%%
\begin{appendix}
\section{The proof of technical lemmas.}\label{app1}
%%%%%%%%%%%%%%%%%%%%%%%%%%%%%%%%%%%%%%%%%%%
\subsection{Lemma \ref{thm:l1}}
%%%%%%%%%%%%%%%%%%%%%%%%%%%%%%%%%%%%%%%%%%%
Positive definiteness: for any $\phi \in L^2(t,\infty)$,
\[
(\phi, T \phi)=\int_0^\infty \frac{dz}{\pi} \left(\int_t^\infty dx e^{-(z+x)^2}\phi(x)\right)^2\geq 0.
\]
Therefore, we have the following estimate for the spectral radius: for any $m=1,2,3,\ldots$,
\bea\label{eq:gel}
\rho_t \leq \left(\tr T^m\right)^{\frac{1}{m}}
\eea
For $t>0$ we can set $m=1$. An explicit calculation gives
\begin{eqnarray}
\rho_t&\leq& \tr T=\frac{1}{\pi}\int_0^\infty dx \int_t^\infty dy e^{-2(x+y)^2}\nonumber\\
&\leq& \frac{1}{\pi}\int_0^\infty dx e^{-2x^2}\int_t^\infty dy e^{-2y^2}
\leq \frac{1}{8}e^{-2t^2}.
 \label{trace}
\end{eqnarray}
Bound (\ref{eq:pbound}) is proved.

For $t<0$, let us re-write $\tr T^m$ in terms 
of the Gaussian random walk:
\begin{eqnarray*}
\tr T^m=\int_{-\infty}^{-t} dy \pr \big[y+B_{2m} \in dy;y+B_{2k+1}\geq 0,\\k=0,1\ldots, m-1;
y+B_{2k}\leq -t, 1\leq k \leq m-1\big].
\end{eqnarray*}
Clearly,
\begin{eqnarray*}
\tr T^m\leq \int_{-\infty}^{-t} dy \pr \big[y+B_{2m} \in dy; y+B_{1}\geq 0
\big].
\end{eqnarray*}
A simple bound on the double integral in the right hand side gives
\[
\tr T^m \leq \frac{1}{\sqrt{2\pi m}} \left(-t+\frac{1}{2\sqrt{2\pi}}e^{-t^2} \right).
\]
Substituting this bound into (\ref{eq:gel}) and optimising with respect to $m$ we arrive
at (\ref{eq:nbound}).
%%%%%%%%%%%%%%%%%%%%%%%%%%
\subsection{Lemma \ref{thm:comb}}
%%%%%%%%%%%%%%%%%%%%%%%%%%%
Let $X_1, X_2, \ldots, X_n$ be the Gaussian increments of a random walk 
$\tilde{B}_k=\sum_{m=1}^k X_m,~k=1,2,\ldots, n$, conditioned to have $\tilde{B}_n=0$. In total, there are $n!$ associated
random walks corresponding to all permutations of the increments. All these permutations are equiprobable. 
Let us say that two random walks are equivalent if one can be obtained from the other by a cyclic
permutation of the increments. Each equivalence class contains $n$ random walks.
Almost surely, precisely one representative in each of the classes will stay negative at all times between $1$
and $n-1$. To construct such a random walk, let us pick an arbitrary representative 
\[
\tilde{B}=(X_1, X_1+X_2, \ldots, X_1+X_2+\ldots+X_{n-1},0)
\]
of the equivalence class. Let $k$ be the point of global maximum of $\tilde{B}$. Then
the following cyclic permutation of $\tilde{B}$ is a random walk which stays below the origin:
\[
(X_{k+1},X_{k+1}+X_{k+2},\ldots, X_{k+1}+X_{k+2}+\ldots+X_n+X_1+\ldots+X_{k-1},0).
\]
All other cyclic permutations of $\tilde{B}$ will have at least one point above the origin.

We conclude that the fraction of random walks conditioned to finish at zero which stay
below the origin at all times between $1$ and $n-1$ is equal to $(n!/n)/n!=1/n$. Therefore,
\[
\pr[\tilde{B}_1<0, \ldots, \tilde{B}_{n-1}<0| \tilde{B}_n =0 ]=\frac{1}{n}.
\]

%%%%%%%%%%%%%%%%%%%%%%%%%%%%%%%%%%%%%
\subsection{Lemma \ref{thm:exit}} 
%%%%%%%%%%%%%%%%%%%%%%%%%%%%%%%%%%%%%
Set $\tau=\tau_0\wedge \tau_t$. By Wald's 
identity \cite{feller} (or the optional stopping theorem for random walks):
\bea
\mathbbm{E}(B_{\tau})=\mathbbm{E}(B_0)=0.
\eea
Therefore,
\[
\mathbbm{E}(B_{\tau}\mathbbm{1}_{\tau_0>\tau_t})+\mathbbm{E}(B_{\tau}\mathbbm{1}_{\tau_0<\tau_t})=0.
\]
Let 
\[
L_{\tau_0}=B_{\tau_0},~ L_{\tau_t}=t-B_{\tau_t}
\]
be the positive overlaps (or `ladder heights') of the random walk  at the exit times $\tau_0$ and $\tau_t$
correspondingly.
It follows from the Wald identity that 
\[
\mathbbm{E}(\mathbbm{1}_{\tau_0>\tau_t})=\frac{1}{|t|}\left(\mathbbm{E}(L_{\tau_0})-\mathbbm{E}((L_{\tau_0}+L_{\tau_t})\mathbbm{1}_{\tau_0>\tau_t})\right)
\]

Applying Spitzer's theorem (see e. g. \cite{feller}, Chapter XVIII.5) to the Gaussian random walk, we
find that
\[
\mathbbm{E}(B_{\tau_0})=\frac{1}{\sqrt{2}},
\]
which implies that
\[
\pr[\tau_0>\tau_t]=\mathbbm{E}(\mathbbm{1}_{\tau_0>\tau_t})=\frac{1}{|t|}\left(\frac{1}{\sqrt{2}}-\mathbbm{E}((L_{\tau_0}+L_{\tau_t})\mathbbm{1}_{\tau_0>\tau_t})\right)
\]
As the second term on the right hand side of the above formula is non-negative, the upper bound
on $\pr[\tau_0>\tau_t]$ claimed in Lemma \ref{thm:exit} is proved.

To establish the lower bound, we apply Cauchy-Schwarz inequality to bound $\mathbbm{E}((L_{\tau_0}+L_{\tau_t})\mathbbm{1}_{\tau_0>\tau_t})$. Solving the resulting inequality with respect to $\mathbbm{E}(\mathbbm{1}_{\tau_0>\tau_t})$ we find that
\[
\mathbbm{E}(\mathbbm{1}_{\tau_0>\tau_t})\geq \frac{1}{\sqrt{2}|t|}\left(1-\sqrt{2\sqrt{2}\mathbbm{E}(L_{\tau_0}^2+L_{\tau_t}^2)}|t|^{-1/2}\right).
\]
The upper bound of Lemma \ref{thm:exit} will be proved with 
$C=\sqrt{2\sqrt{2}\mathbbm{E}(L_{\tau_0}^2+L_{\tau_t}^2)}$ if we
can show that the moments $\mathbbm{E}(L_{\tau_0}^2)$ and $\mathbbm{E}(L_{\tau_t}^2)$
exist. The existence of moments is an immediate 
consequence of the following observation: let $W_s\sim N(0,s)$. Then, for any positive function 
$f$ on $\R$,
\bea\nonumber
\mathbbm{E}(f(W_1)\mathbbm{1}_{W_1>0})\leq  \{\mathbb{E}(f(B_{\tau_0})), ~\mathbb{E}(f(B_{\tau_t})) \}\leq 
2\sup_{s\in (0,1)}\mathbb{E}(f(W_s)\mathbbm{1}_{W_s>0}).\\
\label{eqn:double}
\eea
The desired result follows from the finiteness of the second moment of a Gaussian distribution.

To prove the lower bound in  (\ref{eqn:double}), notice first that
\[
\mathbb{E}(f(B_{\tau_0}))=\mathbb{E}(f(\tilde{B}_{t_0})),
\]
where $(\tilde{B}_n)_{n \geq 0}$ is the Gaussian random walk with $N(0,1)$ increments
and $t_0$ is the exit time through $0$. Then
\[
\mathbb{E}(f(B_{\tau_0}))=\sum_{m=1}^\infty \mathbbm{E}(f(\tilde{B}_m) \mathbbm{1}_{\tau_0=m})
\geq \mathbbm{E}(f(\tilde{B}_1) \mathbbm{1}_{\tau_0=1})=\mathbbm{E}(f(W_1) \mathbbm{1}_{W_1>0}),
\]
which proves the lower bound in (\ref{eqn:double}).

To derive the upper bound, let us consider the standard Brownian motion $(W_s)_{s\geq 0}$
coupled to $(\tilde{B}_n)_{n\geq 0}$ in such a way that $W_s=\tilde{B}_s$ for $s=1,2,\ldots$.
Let $(F^W_s)_{s\geq 0}$ be its natural filtration.
Let
\[
\sigma_{m}=\inf\{s> m-1: W_s=0\}, m=1,2,\ldots
\]
Then applying the strong Markov property,
\begin{eqnarray}
\mathbb{E}(f(B_{\tau_0}))&=&\sum_{m=1}^\infty  \mathbbm{E}(f(\tilde{B}_m) \mathbbm{1}_{\tau_0=m})
=\sum_{m=1}^\infty  \mathbbm{E}((f(\tilde{B}_m)\mathbbm{1}_{\tilde{B}_m>0}) \mathbbm{1}_{\tau_0>m-1})
\nonumber\\
&=&
\sum_{m=1}^\infty  \mathbbm{E}((f(\tilde{B}_m)\mathbbm{1}_{\tilde{B}_m>0}) \mathbbm{1}_{\tau_0>m-1}\mathbbm{1}_{\sigma_m \in (m-1,m)})\nonumber\\
&=&
\sum_{m=1}^\infty  \mathbbm{E}(\mathbbm{E}((f(\tilde{B}_m)\mathbbm{1}_{\tilde{B}_m>0}) \mathbbm{1}_{\tau_0>m-1}\mathbbm{1}_{\sigma_m \in (m-1,m)}\mid F^W_{\sigma_m}))\nonumber\\
&=&
\sum_{m=1}^\infty  \mathbbm{E}(\mathbbm{E}(f(\tilde{B}_m)\mathbbm{1}_{\tilde{B}_m>0})\mid  F^W_{\sigma_m} ) \mathbbm{1}_{\tau_0>m-1}\mathbbm{1}_{\sigma_m \in (m-1,m)})\nonumber\\
&=&
\sum_{m=1}^\infty  \mathbbm{E}(\mathbbm{E}(f(W_{m-\sigma_m})\mathbbm{1}_{W_{m-\sigma_m}>0})\mid  F^W_{\sigma_m} ) \mathbbm{1}_{\tau_0>m-1}\mathbbm{1}_{\sigma_m \in (m-1,m)})\nonumber\\
&\leq &
\sup_{\tau \in (0,1)} \mathbbm{E}(f(W_{\tau})\mathbbm{1}_{W_{\tau}>0})
\sum_{m=1}^\infty  \mathbbm{E}( \mathbbm{1}_{\tau_0>m-1}\mathbbm{1}_{\sigma_m \in (m-1,m)})\label{eq:long}
\end{eqnarray}
By the Reflection Principle,
\[
\mathbbm{E}(\mathbbm{1}_{\tau_0=m}\mathbbm{1}_{\sigma_m\in (m-1,m)})=\frac{1}{2}
\mathbbm{E}(\mathbbm{1}_{\tau_0>m-1}\mathbbm{1}_{\sigma_m\in (m-1,m)}).
\]
Substituting this result into (\ref{eq:long}) and summing over $m$'s we find that
\[
\mathbb{E}(f(B_{\tau_0}))\leq 2 \sup_{\tau \in (0,1)} \mathbbm{E}(f(W_{\tau})\mathbbm{1}_{W_{\tau}>0}).
\]
The proof of the upper bound in (\ref{eqn:double})  for $\mathbbm{E}(f(B_{\tau_0}))$ is complete. 
The derivation of bounds for $\mathbbm{E}(f(B_{\tau_t}))$ is a carbon copy of the above proof.
Lemma \ref{thm:exit} is proved.

%%%%%%%%%%%%%%%%%%%%%%%%%%%%%%%%%%%
%\subsection{Lemma 4}\label{proofpr1}
%%%%%%%%%%%%%%%%%%%%%%%%%%%%%%%%%%%

%%%%%%%%%%%%%%%%%%%%%%%%%%%%%%%%%
\section{A note on the Rider-Sinclair Theorem.}\label{app2}
%%%%%%%%%%%%%%%%%%%%%%%%%%%%%%%%%
For all the definitions, notations and numberings used in this appendix 
we refer readers
to the original paper \cite{rider_sinclair_2014}. Here we deal only with the 
case of even 
sized matrices, therefore we fix the problems happening in paragraph 4.1 
of the paper.

Notice that the only difference between the original 
statement of the Rider-Sinclair theorem and Theorem \ref{thm:rs}
is in the factor $\Gamma_t$. The origin of this 
factor is a 
finite rank perturbation of the operator $T_n\chi$ (see 
\cite[(4.9)]{rider_sinclair_2014}).
In summary, we will correct the following two errors:

(i) the limits of the functions $\tilde{\phi}_n, \tilde{\psi}_n$ (see 
 \cite[(4.15)]{rider_sinclair_2014})
were calculated incorrectly, and this leads to new definitions for $G\left(x\right)$ 
and $g\left(x\right)$ and some 
changes to the constants contained in the third term of the $\Gamma_t$ 
expression as well;

 (ii) it was noted in \cite[p.1644, last paragraph]{rider_sinclair_2014} that
 $\int_{-\infty}^{t} \tilde{\psi}_n \left(x\right)dx$ converges to 
 $G\left(t\right)$ which is true only up to an additive constant and 
 yields corrections to the fourth term of the $\Gamma_t$ expression.
 
 We start with the correct derivation of the limits 
 $\lim\limits_{n\to\infty} \tilde{\phi}_n\left(x\right)$ and 
 $\lim\limits_{n\to\infty} \tilde{\psi}_n\left(x\right)$.
 
	 \begin{eqnarray*}
	 	\tilde{\phi}_n\left(x\right) &=& \kappa_n\int_{0}^{x+\sqrt{n}}
	 	u^{n-2}e^{-u^2/2}d u\\
	 	&=& \kappa_n 2^{\frac{n-3}{2}}\Gamma\left(\frac{n-1}{2}\right)
	 	P\left(\frac{n-1}{2},\frac{\left(x+\sqrt{n}\right)^2}{2}\right),
	 \end{eqnarray*}
 where 
 $P\left(a,x\right)=\frac{\Gamma\left(a,x\right)}{\Gamma\left(a\right)}$
 is the incomplete regularized Gamma function. It is easy to check by using
 the duplication formula for Gamma functions that the prefactor in front of $P$ 
 is asymptotically equal to $2^{-1/2}$. For the last factor we 
 use the well known asymptotic formula (see, e.g. 
 \cite[(9.16)]{borodin_sinclair})
	 \begin{equation*}
	 	P\left(a,a+\sqrt{2a}x\right) \sim \frac{1}{2}\erfc\left(-x\right), 
	 	a\to\infty.
	 \end{equation*}
 Together with the above identity we get
	 \begin{equation}\label{eq:G_def}
	 	G\left(x\right) = \lim\limits_{n\to\infty} 
	 	\tilde{\phi}_n\left(x\right)
	 	=\frac{1}{2\sqrt{2}}\erfc\left(-x\right) = 
	 	\int_{-\infty}^{x}\frac{e^{-t^2}}{\sqrt{2\pi}} d t.
	 \end{equation}
 For $\tilde{\psi}_n$ we have
	 \begin{eqnarray*}
	 	\tilde{\psi}_n\left(x\right) &=& 
	 	\kappa_n'\left(\sqrt{n}+x\right)^{n-1}
	 	e^{-\left(x+\sqrt{n}\right)^2/2} \\&=&
	 	\kappa_n'n^{\frac{n-1}{2}}e^{-n/2}e^{\left(n-1\right)
	 	\log\left(1+\frac{x}{\sqrt{n}}\right)-x\sqrt{n}-x^2/2}\\
	 	&=&
	 	\kappa_n'n^{\frac{n-1}{2}}e^{-n/2}
		e^{-x^2+O\left(n^{-1/2}\right)}.
	 \end{eqnarray*}
The prefactor can be shown asymptotically equal to $\frac{1}{\sqrt{2\pi}}$
 and therefore
	 \begin{equation}\label{eq:g_def}
	 	g\left(x\right)=\lim\limits_{n\to\infty} 
	 	\tilde{\psi}_n\left(x\right) = \frac{1}{\sqrt{2\pi}}e^{-x^2}.
	 \end{equation}
 \begin{remark}
 	Notice that the the correct expression for the product
	 	\begin{equation*}
	 		g\left(x\right)G\left(y\right) = \dfrac{1}{4\sqrt{\pi}}
	 		e^{-x^2}\erfc\left(-y\right)
	 	\end{equation*}
	 is consistent with the scaling limit of the kernel 	 
	 $S_n\left(x,y\right)$ 
	 from \cite{borodin_sinclair}.
 \end{remark}
 
  \begin{remark}
  	Throughout this appendix we use the definitions 
  	\eqref{eq:G_def}, \eqref{eq:g_def} for the functions $g$ and $G$ which differ from the conventions adopted
  	in Theorem~\ref{thm:rs}. This is done to highlight the changes
  	needed to calculate the correct scaling limit of 
  	the $2\times 2$ determinant 
  	$\det\left(1-\left(\alpha_i,\beta_j\right)\right)_{1\leq i,j \leq 2}$  	
  	(see \cite[p.1641-1644]{rider_sinclair_2014}). For convenience, we placed the terms correcting original answers in boxes.
  \end{remark}

 The first matrix element $\left(\alpha_1,\beta_1\right)$ converges to
 	\begin{equation*}
 		\int_t^\infty G(x) (I-T\chi_t)^{-1} g(x)dx,
 	\end{equation*}
 with operators $T$ and $\chi_t$ defined as in Theorem~\ref{thm:rs}.
 The second term $\left(\alpha_1,\beta_2\right)$ converges to
 	\begin{equation*}
 		(I-T\chi_t)^{-1}G(t) - G\left(\infty\right) = 
 		(I-T\chi_t)^{-1}G(t) - \boxed{\frac{1}{\sqrt{2}}}.
 	\end{equation*}
To calculate the remaining two terms (see \cite[(4.12)]{rider_sinclair_2014}) we 
follow 
\cite[p. 1644, last paragraph]{rider_sinclair_2014} to arrive at
	\begin{eqnarray*}
		\left(\alpha_2,\beta_1\right) &=& \frac{1}{2}\int_{-\infty}^{t}
		\left(1-\tilde{T}_n\chi_t\right)^{-1}
		\tilde{\psi}_n\left(x\right)d x\\
		&=& \frac{1}{2}\int_{-\infty}^{t}
		\tilde{T}_n\chi_t\left(1-\chi_t\tilde{T}_n\chi_t\right)^{-1}
		\tilde{\psi}_n\left(x\right)d x
		+ \frac{1}{2}\int_{-\infty}^{t}\tilde{\psi}_n\left(x\right)d x.
	\end{eqnarray*}
The first summand on the right hand side converges to its formal limit
	\begin{equation*}
		\frac{1}{2}\int_{-\infty}^{t}
		T\chi\left(1-\chi_t T_n\chi_t\right)^{-1}
		g\left(x\right)d x.
	\end{equation*}
The answer for the second summand depends on the parity of $n$:
	\begin{equation*}
		\int_{-\infty}^{t}\tilde{\psi}_n\left(x\right)d x 
		\approx -\frac{1}{2\sqrt{2}} \erfc\left(t\right)
		+ \frac{1-\left(-1\right)^n}{\sqrt{2}} = 
		G\left(t\right) - \frac{\left(-1\right)^n}{\sqrt{2}}.
	\end{equation*}
Considering only the even sized matrices we conclude that 
$\left(\alpha_2,\beta_1\right)$ converges to
	\begin{equation*}
		\frac{1}{2}\int_{-\infty}^t (I-T\chi_t)^{-1} g(x)dx-\boxed{\frac{1}{2\sqrt{2}}}.
	\end{equation*}
Therefore, the limit of $\left(\alpha_2,\beta_2\right)$  
coincides with the answer stated in the original paper and is equal to
	\begin{equation*}
		\frac{1}{2}\int_{-\infty}^t R(x,t+)dx.
	\end{equation*}
Gathering the answers derived above we find that
	\begin{eqnarray*}
		\Gamma_t &=& \left(1-\int_t^\infty G(x) (I-T\chi_t)^{-1} g(x)dx\right)
		\left(1-\frac{1}{2}\int_{-\infty}^t R(x,t+)dx\right)
		\\
		&+&\frac{1}{2}
		\left(
		\frac{1}{\sqrt{2}} -\left(I-T\chi_t\right)^{-1}G(t)
		\right)
		\left(
		\int_{-\infty}^t (I-T\chi_t)^{-1} g(x)dx-\frac{1}{\sqrt{2}}
		\right),
	\end{eqnarray*}
which is exactly the statement of Theorem~\ref{thm:rs} after the re-definition
$g,G \to \frac{1}{\sqrt{2}}g,\frac{1}{\sqrt{2}}G$.
\end{appendix}
%%%%%%%%%%%%%%%%%%%%%%%%%%%%%%%%

\end{document}